\documentclass[a4paper]{article}

\usepackage[utf8]{inputenc} 
\usepackage[T1]{fontenc} 
\usepackage{geometry} 
\usepackage[english]{babel} 

\title{Sphere paths in outer space}
\author{Camille Horbez}
\date{}

\usepackage{graphicx}
\usepackage[all]{xy}
\usepackage{amsmath,amsfonts, amssymb,amsthm}
\usepackage{bbm}
\begin{document}
\large
\maketitle
\newtheorem{de}{Definition} [section]
\newtheorem{theo}[de]{Theorem} 
\newtheorem{prop}[de]{Proposition}
\newtheorem{lemma}[de]{Lemma}
\newtheorem{cor}[de]{Corollary}
\newtheorem{propd}[de]{Proposition-Definition}

\begin{abstract}    
We give estimates on the length of paths defined in the sphere model of outer space using a surgery process, and show that they make definite progress in some sense when they remain in some thick part of outer space. To do so, we relate the Lipschitz metric on outer space to a notion of intersection numbers.
\\
\\
\end{abstract}

\maketitle

\section*{Introduction}

In order to study the outer automorphism group of a finitely generated free group, Culler and Vogtmann introduced a space, called \emph{outer space}, on which the group $Out(F_n)$ acts in a nice way (\cite{CV}, see also \cite{Vog} for a good survey article). This space is built as an analog of Teichmüller spaces, used to study the mapping class group of a surface. While Teichmüller spaces are equipped with several interesting metrics, whose properties have been investigated a lot, there had been no systematic investigation of metric properties of outer space before Francaviglia and Martino studied an analog of Thurston's asymmetric metric (\cite{FM}). In particular, Francaviglia and Martino proved that outer space is geodesic for this metric, the geodesics being obtained by using a folding process. 

Building on ideas of Whitehead (\cite{Whi}), Hatcher defined a new model for outer space, using sphere systems in a $3$-dimensional manifold with fundamental group $F_n$ (\cite{Hat}). In order to prove the contractibility of the full sphere complex, he also defined a combing path in this model of outer space, which appears to look like an "unfolding path". A modification of this path was also used by Hatcher and Vogtmann to prove exponential isoperimetric inequalities for $Out(F_n)$ (\cite{HV}).

Our goal is to investigate the metric properties of this path. As combing paths are piecewise linear, we can talk about their vertices, and define the \emph{length} $l(\gamma)$ of a combing path $\gamma$ to be the sum of the distances from one vertex to the next. We prove the following result.
\\
\\
\textbf{Main theorem} : \textit{For all $n\ge 2$ and $\epsilon>0$, there exist $K,L\in\mathbb{R}$ such that the following holds.}

\noindent \textit{Let $A,B\in CV_n$ be such that the combing path $\gamma$ from $B$ to $A$ remains in the $\epsilon$-thick part of outer space. Then}

\begin{center}
\textit{$\frac{l(\gamma)}{K}-L\le d(A,B)\le l(\gamma).$}
\end{center}

\indent In section \ref{sec1}, we recall some basic facts about the different models of outer space. In particular, we recall the definition of the Lipschitz metric on outer space (section \ref{sec1.1}), as well as two notions of intersection numbers on outer space : the first was introduced by Guirardel in \cite{Gui}, who defined a convex core for two actions of groups on trees (section \ref{sec1.2}), and the second is a geometric notion of intersection in the sphere model of outer space (section \ref{sec1.3}). Intersection numbers have turned out to be a powerful tool in the study of mapping class groups - they were used for example by Bowditch to give a new proof of the hyperbolicity of the curve complex (\cite{Bow}). It seems that they are also relevant to study paths in the sphere model of outer space. Finally, in section \ref{sec1.4}, we recall the definition of the combing path in the sphere model of outer space from the works of Hatcher (\cite{Hat}) and Hatcher-Vogtmann (\cite{HV}).

In section \ref{sec2}, we establish the equality between both notions of intersection numbers (section \ref{sec2.1}). We then prove the following relation between intersection numbers and the Lipschitz metric (with the convention $\log 0=0$), which may be of independent interest. 
\\
\\
\noindent \textbf{Theorem} :
\textit{For all $n\ge 2$ and $\epsilon>0$, there exist $K',L'\in\mathbb{R}$ such that for all points $X,Y$ in the $\epsilon$-thick part of $CV_n$, we have
\begin{displaymath}
\frac{1}{K'}\log(i(X,Y))-L'\le d(X,Y)\le K' \log(i(X,Y))+L'.
\end{displaymath}}

\indent Section \ref{sec3} is dedicated to the proof of our main theorem. The main step in our proof is to understand the growth of intersection numbers along combing paths (section \ref{sec3.2}). We prove the following estimate, which can be regarded as an analog for combing paths of the result of Behrstock, Bestvina and Clay about growth of intersection numbers along the axis of a fully irreducible automorphism of $F_n$ (\cite{BBC}). 
\\
\\
\noindent \textbf{Proposition} : \textit{For all $n\ge 2$ and $\epsilon>0$, there exist $C_1, C_2\in (1,+\infty)$ such that the following holds.
\\ Let $A,B\in CV_n$ be such that the combing path from $B$ to $A$ stays in the $\epsilon$-thick part of $CV_n$, and let $A=A_0,\dots,A_N=B$ be the vertices of this path. Assume that $N\ge 3$, then }

\begin{displaymath} 
C_1^{N}\le i(A,B)\le C_2^{N}.
\end{displaymath}

\noindent The main theorem easily follows from the two results above.
\\
\\
\textbf{Acknowledgments} : It is a pleasure to thank Karen Vogtmann, who gave me the opportunity to come to Cornell University to work with her, and took a lot of her time to introduce me to this wonderful area of mathematics. None of this work would have been possible without her valuable help and comments.

\section{Preliminaries} \label{sec1}

\subsection{Outer space in terms of graphs and the Lipschitz metric} \label{sec1.1}

We recall the construction of outer space by Culler and Vogtmann (\cite{CV}). 

A \emph{metric graph} is a graph, all of whose vertices have valence at least three, endowed with a path metric : each edge is assigned a positive length $l$ that makes it isometric to the segment $[0,l]$ in $\mathbb{R}$. Denote by $R_n$ the metric graph (called a \emph{rose}) with one vertex and $n$ edges of length $\frac{1}{n}$. A \emph{marking} of a metric graph $G$ of fundamental group $F_n$ is a homotopy equivalence $\rho:R_n\to G$. Define an equivalence relation on the collection of marked metric graphs by $(G,\rho)\sim (H,\rho')$ if there exists a homothety $h:G\to H$ such that $h\circ\rho$ is homotopic to $\rho'$. \emph{Outer space}, denoted by $CV_n$, is defined to be the set of classes of marked metric graphs under this equivalence relation. As we took the quotient by homotheties, we can assume the graphs to be normalized to have total length $1$.

To every marked metric graph $G$, one associates an open simplex by making the lengths of the edges of $G$ vary, with sum equal to $1$. The simplex of a graph $H$ is identified with a face of the simplex of a graph $G$ if $H$ can be obtained from $G$ by shrinking the lengths of some edges to $0$. Outer space is endowed with the quotient topology of the natural topology on the union of the simplices by these face identifications.

Given $\epsilon>0$, the \emph{$\epsilon$-thick part} of outer space is the subspace consisting of graphs (normalized to have length $1$) that do not contain a loop of length less than $\epsilon$.

The group $Out(F_n)$ acts on outer space by precomposing the markings. More precisely, let $(G,\rho)$ be a marked metric graph and $\Phi\in Out(F_n)$. Choose some representative $\phi\in Aut(F_n)$ for $\Phi$, and a homotopy equivalence $f:R_n\to R_n$ that induces $\phi$ on the fundamental group. The action of $\Phi$ on $(G,\rho)$ is given by $[G,\rho]\Phi=[G,\rho\circ f]$ (it is easy to check that this definition does not depend on the choices of $\phi$ and $f$). This action is not cocompact. However, the group $Out(F_n)$ acts cocompactly on the \emph{spine} of outer space, which is defined to be the geometric realization of the poset of the simplices of outer space, ordered by inclusion of faces.
\\
\\ 
\indent In \cite{FM}, Francaviglia and Martino defined a metric on outer space, compatible with the topology defined above, in the following way. Given two marked metric graphs $(G,\rho)$ and $(H,\rho')$, a \emph{difference of markings} from $G$ to $H$ is a map which is homotopic to $\rho'\circ\rho^{-1}$, where $\rho^{-1}$ denotes a homotopy inverse of the homotopy equivalence $\rho$. Define the \emph{stretching factor} $\Lambda(G,H)$ from $G$ to $H$ as the infimum of the Lipschitz constant of a difference of markings from $G$ to $H$. Francaviglia and Martino proved (\cite[Theorem 4.17]{FM}) that $d(G,H):=\log \Lambda(G,H)$ defines an asymmetric metric on outer space, which is $Out(F_n)$-invariant.
This metric is not symmetric, and not even quasi-symmetric (see e.g. the examples in section 1.3 of \cite{AKB}). However, it is quasi-symmetric when restricted to the $\epsilon$-thick part of outer space for some $\epsilon>0$. 

\begin{prop} (\cite[Proposition 1.3]{HM}; \cite[Theorem 24]{AKB}) \label{prop1.1}
For all $\epsilon>0$, there exists $C=C(\epsilon)>0$ such that for all $X,Y$ in the $\epsilon$-thick part of outer space, we have $d(Y,X)\le Cd(X,Y)$, i.e. $\Lambda(Y,X)\le \Lambda(X,Y)^C$.
\qed
\end{prop}

\subsection{Actions on trees and Guirardel's intersection number} \label{sec1.2}

The universal cover of a marked metric graph is a metric tree, endowed with an action of $F_n$ given by the marking. This action is free and isometric. It is also \emph{minimal}, meaning that there is no proper invariant subtree.  One can define outer space as the set of all minimal, free, isometric actions of $F_n$ on metric simplicial trees, up to equivariant homothety.
\\
\\
\indent In \cite{Gui}, Guirardel defined a notion of intersection number between two actions on trees. We recall his construction. Let $T_1,T_2$ be two simplicial metric trees with free, minimal, isometric actions of $F_n$. A \emph{direction} in $T_1$ is a component of $T_1-\{x\}$, for some point $x\in T_1$. A \emph{quadrant} in $T_1\times T_2$ is the product $\delta_1\times\delta_2$ of a direction $\delta_1$ in $T_1$ and a direction $\delta_2$ in $T_2$. 
\\
\\
\indent Let $*_1$ (resp. $*_2$) be a fixed basepoint in $T_1$ (resp. $T_2$). A quadrant $Q=\delta_1\times\delta_2$ is said to be \emph{heavy} if there exists a sequence $(g_k)$ of elements in $F_n$ such that :

1) $g_k(*_1,*_2)\in Q$,

2) $\lim_{k\to +\infty} d_{T_1}(*_1,g_k*_1)=+\infty$ and $\lim_{k\to +\infty}d_{T_2}(*_2,g_k*_2)=+\infty$.

\noindent Otherwise $Q$ is said to be \emph{light}.
\\
\\
\noindent \textit{Remarks} : 1) This definition does not depend on the choice of the basepoint in $T_1\times T_2$. 
\\
\\
2) For every $g\in F_n$, a quadrant $Q$ is heavy if and only if its translate $gQ$ is heavy.
\\
\\
\indent The \emph{core} $\mathcal{C}(T_1\times T_2)$ of $T_1\times T_2$ is defined to be the complement of the union of all light quadrants in $T_1\times T_2$. By the second remark above, it is an $F_n$-invariant subset of $T_1\times T_2$. The \emph{intersection number} $i(T_1,T_2)$ between $T_1$ and $T_2$ is the number of $2$-cells in $\mathcal{C}(T_1\times T_2)/F_n$. (This definition is slightly different from Guirardel's, as it does not take into account the lengths of the edges of the trees. In other words, we consider $T_1$ and $T_2$ as simplicial trees with all edges having length $1$ for computing the intersection number).
 Given an edge $e_1\subset T_1$, the \emph{slice} of the core at $e_1$ is 
 
\begin{center} 
 $\mathcal{C}_{e_1}=\{e_2\in T_2|e_1\times e_2\subset\mathcal{C}(T_1\times T_2)\}$.
\end{center}

\noindent Let $e_1\subset T_1$ be an edge, and $g\in F_n$. By $F_n$-invariance of the core, we have $\mathcal{C}_{ge_1}=g\mathcal{C}_{e_1}$. The intersection number is thus equal to

 \begin{displaymath}
i(T_1,T_2)=\sum_{e_1\subset T_1/F_k}|\mathcal{C}_{e_1}|,
\end{displaymath}

\noindent where $|\mathcal{C}_{e_1}|$ denotes the cardinality of $\mathcal{C}_{e_1}$.
\\
\\
\indent In \cite[section 3]{BBC}, Behrstock, Bestvina and Clay gave an algorithm to compute the slices of the core $\mathcal{C}(T_1,T_2)$ for $T_1,T_2\in CV_n$, and hence the intersection number $i(T_1,T_2)$. We now describe their construction, which we will use in the proof of proposition \ref{prop2.8} to compare the intersection number with the metric on outer space.

Let $f:T_1\to T_2$ be a \emph{morphism}, i.e. an equivariant cellular map that linearly expands each edge in $T_1$ over a tight edge path in $T_2$ (note that this definition of a morphism between trees is slightly different from the usual one). It descends to a homotopy equivalence $\sigma:\Gamma_1\to\Gamma_2$, where $\Gamma_1$ (resp. $\Gamma_2$) is the graph corresponding to $T_1$ (resp. $T_2$) in outer space, i.e. its quotient by $F_n$. Fix a morphism $f':T_2\to T_1$ such that $\sigma':\Gamma_2\to\Gamma_1$ is a homotopy inverse of $\sigma$. 
Fix basepoints $*_1\in T_1$ and $*_2\in T_2$ such that $f'(*_2)=*_1$. Slightly abusing notations, we will again denote by $*_1$ and $*_2$ their projections to $\Gamma_1$ and $\Gamma_2$.

Let $e$ be an oriented edge of $\Gamma_1$. Subdivide $e$ into $e_+e_-$, and let $p_e$ be the subdivision point. 
Fix a tight edge path $\alpha_e\subset\Gamma_1$ from $*_1$ to $p_e$ which has final edge $e_+$. Let $\Sigma_e=(\sigma')^{-1}(p_e)\subset\Gamma_2$. For $q\in\Sigma_e$, there is a tight path $\gamma_q$ in $\Gamma_2$ from $*_2$ to $q$ such that up to homotopy, we have $\alpha_e=[\sigma'(\gamma_q)]$, where $[\sigma'(\gamma_q)]$ denotes the path obtained after tightening $\sigma'(\gamma_q)$. As $\sigma'$ is a homotopy equivalence, the path $\gamma_q$ is unique. 
Let $\widetilde{\gamma_q}$ be the lift of $\gamma_q$ to $T_2$ that originates at $*_2$, let $\widetilde{\Sigma_e}$ be the set of all terminal points of $\widetilde{\gamma_q}$ for $q$ varying in $\Sigma_e$, and let $T_e$ be the subtree of $T_2$ spanned by $\widetilde{\Sigma_e}$. Behrstock, Bestvina and Clay proved the following result (in fact, they even gave an algorithm that enables us to compute precisely the slice of the core at $e$ from the tree $T_e$).

\begin{prop} (\cite[Lemma 3.7]{BBC}) \label{prop1.2}
The slice of the core at $e$ is contained in $T_e$.
\qed
\end{prop}

\subsection{Outer space in terms of sphere systems} \label{sec1.3}

Outer space has a description in terms of sphere systems in a $3$-dimensional manifold with fundamental group $F_n$, which was introduced by Hatcher (\cite{Hat}).

Let $n\in\mathbb{N}$, and $M_n=\#_nS^1\times S^2$ be the connected sum of $n$ copies of $S^1\times S^2$. The fundamental group of $M_n$ is a free group of rank $n$. A \emph{sphere set} is a collection of disjoint embedded $2$-spheres in $M_n$. A \emph{sphere system} $S$ is a sphere set such that no sphere in $S$ bounds a ball in $M_n$, and no two spheres in $S$ are isotopic. A \emph{weighted sphere system} is a sphere system in which each sphere is assigned a positive weight, with the sum of all weights equal to $1$. A sphere set $S$ is said to be \emph{simple} if every component of $M_n-S$ is simply connected. \emph{Outer space} is defined to be the set of all isotopy classes of weighted simple sphere systems. 

The equivalence with the definitions in the previous sections was shown by Hatcher in \cite[Appendix]{Hat}. A simple sphere system $S$ has a dual graph $G(S)$ whose vertices are the components of the complement of $S$ in $M_n$, and whose edges are the spheres in $S$. The graph $G(S)$ can be embedded in $M_n$, each vertex lying in one component of $M_n-S$, and each edge crossing exactly one sphere of $S$ exactly once.
\\
\\
\indent An important tool in the study of sphere systems is Hatcher's normal form. Let $\Sigma$ be a simple sphere system. A sphere system $S$ is said to be in \emph{normal form} with respect to $\Sigma$ if every sphere in $S$ either

\begin{enumerate}
\item belongs to $\Sigma$, or

\item is disjoint from $\Sigma$ and not isotopic to any sphere in $\Sigma$, or 

\item intersects $\Sigma$ transversely in a collection of circles that split it into components called \emph{pieces}, in such a way that for each component $P$ of $M_n-\Sigma$,

(i) each piece in $P$ meets each component of $\partial P$ in at most one circle.

(ii) no piece in $P$ is a disk which is isotopic, fixing its boundary, to a disk in $\partial P$.
\end{enumerate}

The following proposition was first proved by Hatcher when $\Sigma$ is a maximal sphere system (\cite[Propositions 1.1 and 1.2]{Hat}). The extension to the general case is easy, and can be found for example in \cite[Propositions 2.1 and 2.2]{HV}.

\begin{prop} \label{prop1.3}
Every sphere system $S$ is isotopic to a sphere system in normal form with respect to $\Sigma$. Besides, the number of intersection circles between a sphere system $S'$ isotopic to $S$ and $\Sigma$ is minimized if and only if $S'$ is in normal form with respect to $\Sigma$.
\hfill $\square$
\end{prop}

Given $X,Y\in CV_n$, we define their \emph{geometric intersection number} $i(X,Y)$ as the minimal number of intersection circles between a sphere system $S$ representing $X$ and a sphere system $S'$ representing $Y$. This is equal to the number of intersection circles between two representatives in normal form. In the same way, if $s\in X$ and $s'\in Y$ are two spheres, then we define $i(s,s')$ to be the minimal number of intersection circles between a sphere isotopic to $s$ and a sphere isotopic to $s'$. Again, this is achieved when $X$ and $Y$ are in normal form. 

\subsection{Sphere paths in outer space} \label{sec1.4}

Let $\Sigma, S$ be two simple sphere systems in $M_n$, and assume that $S$ is in normal form with respect to $\Sigma$. Following \cite{HV}, we describe a surgery process for producing a new simple sphere system $S'$ from $S$ which intersects $\Sigma$ in fewer circles. Let $C$ be a circle component of $S\cap\Sigma$ which bounds an innermost disk $D\subset \Sigma$, and let $s$ be the sphere of $S$ that contains $C$. Taking a parallel copy of $s$ and performing surgery on it along $D$ creates two spheres $s'$ and $s''$. The new sphere system $S'$ is obtained from $S$ by deleting $s$, replacing it by $s'\cup s''$, and, if necessary, identifying parallel spheres and deleting spheres that bound a ball. We say that the sphere system $S'$ is obtained by performing surgery on $S$ along $\Sigma$. Hatcher and Vogtmann proved (\cite[Lemma 3.1]{HV}) that the sphere system $S'$ is again simple.
\\
\\
\indent Using this surgery process, Hatcher and Vogtmann defined a canonical path between two points in outer space. The idea is to perform simultaneously all surgeries on $S$ along innermost intersection circles in $\Sigma$. However, we have to be careful while defining these processes. Indeed, problems occur when two of the surgery disks lie on different sides of a sphere $s\in S$ (because in that case, it is impossible to choose the parallel copy of $s$ on which we perform surgery), or when one sphere $\sigma\in\Sigma$ intersects $S$ only once (because there are two possible choices for the disk $D$, and we want the construction to be canonical). To solve these problems, Hatcher and Vogtmann use a doubling trick. 

\begin{figure} 
\begin{center}
\includegraphics[width=125mm,height=35mm]{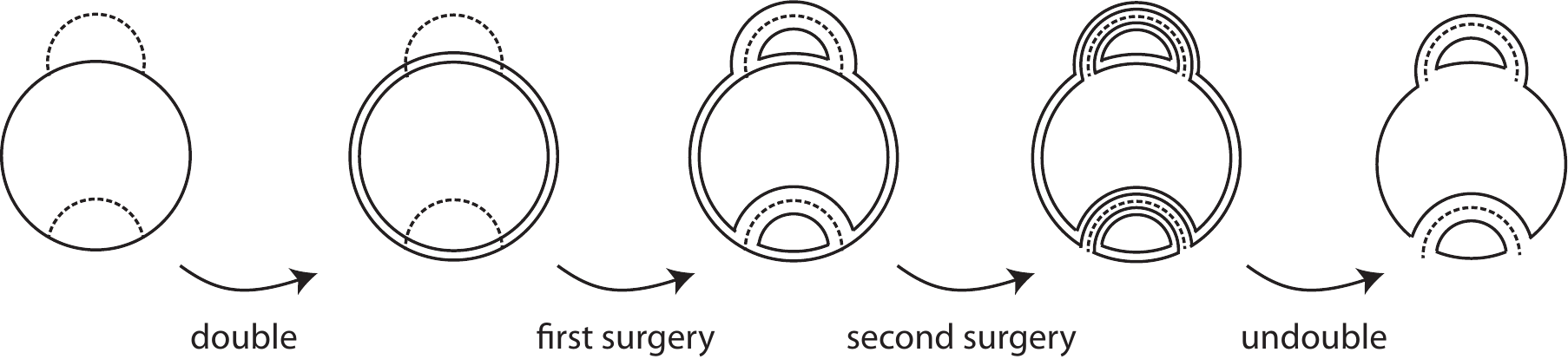}
\caption{A double surgery step}
\label{fig1}
\end{center}
\end{figure}

\begin{figure} 
\begin{center}
\includegraphics[width=10cm,height=4cm]{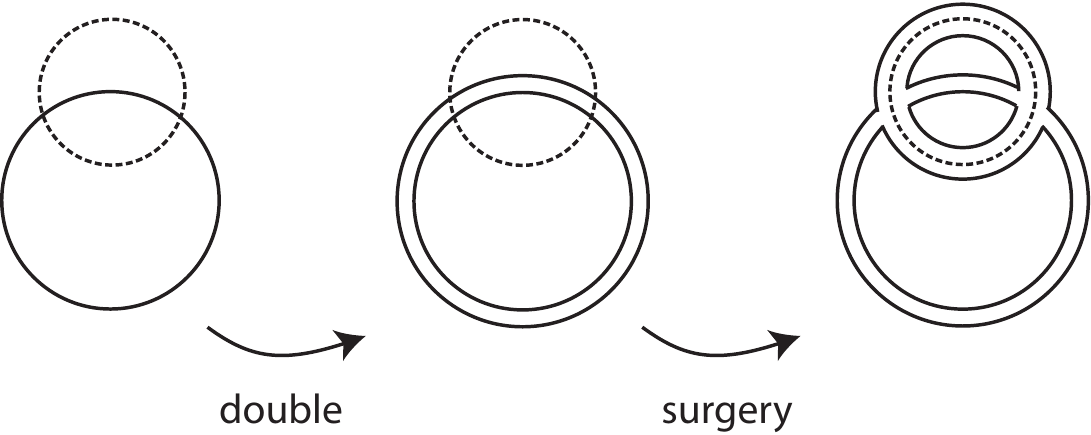}
\caption{Surgery along the last intersection circle with a sphere in $\Sigma$}
\label{fig2}
\end{center}
\end{figure}

Start by adding a parallel copy of each sphere $s\in S$ to get a sphere set $\hat{S}$ (step 1 in Figure \ref{fig1}), and give to each copy of $s$ half of the weight of $s$. Then perform simultaneous surgeries on $\hat{S}$ along all disks in $\Sigma$ that are innermost among the disks bounded by an intersection circle between $\hat{S}$ and $\Sigma$ (step 2 in Figure \ref{fig1}). This operation is now well-defined because on each copy of the sphere, all surgeries are performed on the same side. Besides, all intersection circles have been doubled, so no sphere $\sigma\in\Sigma$ intersects $S$ exactly once. In that way, we get a new sphere set $\hat{S'}$, whose projection $S'$ to outer space after deleting trivial spheres and identifying parallel spheres does not intersect $S$, hence $S$ and $S'$ share a closed simplex in outer space. During the process, transfer continuously the weight of any sphere in $\hat{S}$ on which surgery is performed equally between the nontrivial spheres in $S'$ obtained from it after the surgery. Then perform again simultaneous surgeries on $\hat{S'}$ along all disks in $\Sigma$ that are innermost among the disks bounded by an intersection circle between $\hat{S'}$ and $\Sigma$ to get a sphere set $\hat{S''}$ (step 3 in Figure \ref{fig1}). Again, we transfer continuously the weights of a sphere $s'\in S'$ equally between the spheres that come from it. Assume that no sphere $\sigma\in\Sigma$ intersects $S$ exactly once. Then the sphere set $\hat{S''}$ is the double of a simple sphere system $S''$ : indeed, performing two successive surgery steps on a sphere $s\in B$ consists of performing the "same" surgery on each side of the sphere $s$. We can thus "undouble" the sphere set $\hat{S''}$ (step 4 in Figure \ref{fig1}). In the case when one of the sphere $\sigma\in\Sigma$ intersects $S$ exactly once, it is no longer true that $\hat{S''}$ is a doubled sphere system (see Figure \ref{fig2}). In that case, we can still define $S''$ to be the sphere system we get from $\hat{S''}$ by deleting trivial spheres and identifying parallel spheres if necessary. We repeat this whole process, starting from the sphere system $S''$, until we get a sphere system $\Sigma'$ which does not intersect $\Sigma$ (a slight variation on the argument in \cite{HV} shows that this process eventually stops). We then join $\Sigma'$ to $\Sigma$ by the unique straight line between them in the closed simplex they share. The path we get, which we call the \emph{combing path} from $S$ to $\Sigma$, is a piecewise linear path, which we parametrize by arc length. In the sequel, we will say that the sphere system $S''$ we get by doubling $S$, performing two successive surgeries and undoubling the result is obtained by performing a double surgery step on $S$ along $\Sigma$.

\section{Intersection numbers and the Lipschitz metric} \label{sec2}

\subsection{Equality between both notions of intersection numbers} \label{sec2.1}

In \cite{Gui}, Guirardel proved the equivalence of his notion of intersection number with an algebraic notion of intersection number defined by Scott in \cite{Sco98} in the case of a one-edge splitting of a group. Gadgil and Pandit (\cite{Gad06},\cite{GP09}) then showed that this was also equivalent to the geometric intersection number between two spheres in $M_n$, each sphere corresponding to a one-edge free splitting of $F_n$. In this section, we give a direct proof of the equivalence of the geometric notion of intersection number between two points in outer space and Guirardel's notion. We thank Vincent Guirardel for his help for the proof of this result, which follows the ideas of (\cite[Théorème 7.1]{Gui}).

\begin{prop} \label{prop2.1}
The geometric intersection number coincides with Guirardel's intersection number.
\end{prop}

\noindent \textit{Proof} : Let $S_1$, $S_2$ be two simple sphere systems in normal form. Let $\pi:\widetilde{M_n}\to M_n$ be a universal cover of $M_n$. For $i\in\{1,2\}$, let $\widetilde{S_i}:=\pi^{-1}(S_i)$, and let $T_i$ be the tree dual to $\widetilde{S_i}$, i.e. the tree which has a vertex for each component of $\widetilde{M_n}-\widetilde{S_i}$ and an edge for each component of $\widetilde{S_i}$. Let $A_i$ be a small open tubular neighborhood of $S_i$ (of the form $S_i\times [-\epsilon,\epsilon]$ for some $\epsilon>0$), and let $\widetilde{A_i}$ be its lift to $\widetilde{M_n}$. There is an $F_n$-equivariant map $f_i:\widetilde{M_n}\to T_i$ that maps a component of $\widetilde{A_i}$ to the corresponding edge in $T_i$, and a component of $\widetilde{M_n}-\widetilde{A_i}$ to the corresponding vertex. Let $F=(f_1,f_2):\widetilde{M_n}\to T_1\times T_2$. The preimage of any point $x_1\in T_1$ (resp. $x_2\in T_2$) by $f_1$ (resp. $f_2$) is connected. In addition, it is easy to check that the image $F(\widetilde{M_n})$ is closed. Therefore, by \cite[Corollaire 5.3]{Gui}, we have $\mathcal{C}(T_1\times T_2)\subset F(\widetilde{M_n})$. Furthermore, the map $F$ sends each connected component of $(\widetilde{S_1}\cap\widetilde{S_2})\times [-\epsilon,\epsilon]^2$ to a $2$-cell in $T_1\times T_2$, and conversely the preimage of a $2$-cell in $T_1\times T_2$ is of this form, if it is nonempty. Therefore, we just have to check that intersections between two spheres are mapped by $F$ to cells in the core. 

Let $\Sigma_1$ be a sphere in $\widetilde{S}_1$, and $\Sigma_2$ be a sphere in $\widetilde{S}_2$, such that $\Sigma_1\cap\Sigma_2\ne\emptyset$. Using the fact that $S_1$ and $S_2$ are in normal form, we will show that all the components of the complement of $\Sigma_1\cup\Sigma_2$ in $\widetilde{M_n}$ are unbounded. Suppose by contradiction that there exists a bounded component $X$ of $\widetilde{M_n}-(\Sigma_1\cup\Sigma_2)$. Then the set of components of $\widetilde{M_n}-\widetilde{S_1}$ that intersect $X$ is finite, hence its projection to $T_1$ (by $f_1$) is a finite collection of vertices in $T_1$, one of whose, which we denote by $v_{max}$, must be extremal in the subtree they span. Let $Y$ be the component of $\widetilde{M_n}-\widetilde{S_1}$ which projects to $v_{max}$ by $f_1$. The boundary of $X\cap Y$ intersects only one of the boundary spheres of $Y$. As $S_1$ and $S_2$ are in normal form, the boundary of $X\cap Y$ is a disk whose boundary circle lies in one boundary sphere of $Y$, and which is not isotopic, fixing its boundary, to a disk in $\partial Y$. This implies that $X$ is unbounded.

As a result, for every connected component of $(\widetilde{S_1}\cap\widetilde{S_2})\times [-\epsilon,\epsilon]^2$, and for every quadrant containing the correponding $2$-cell $e_1\times e_2\subset T_1\times T_2$, there exists an element $g\in F_n$ whose iterates $g^k$ map a basepoint in $T_1\times T_2$ towards infinity within the quadrant. So every quadrant containing $e_1\times e_2$ is heavy, hence $e_1\times e_2\subset \mathcal{C}(T_1\times T_2)$.
\qed

\subsection{Intersection numbers and the Lipschitz metric} \label{sec2.2}

In this section, we relate intersection numbers to the Lipschitz metric on outer space. We will show the following asymptotic estimate comparing the intersection number and the distance between two points in the $\epsilon$-thick part of outer space (note that no comparison can hold on the entire outer space, as intersection numbers are symmetric whereas the Lipschitz metric is not). In the following statement, we take the convention that $\log 0=0$.

\begin{theo} \label{theo2.2}
For all $n\ge 2$ and $\epsilon>0$, there exist $K',L'\in\mathbb{R}$ such that for all points $X,Y$ in the $\epsilon$-thick part of $CV_n$, we have
\begin{displaymath}
\frac{1}{K'}\log(i(X,Y))-L'\le d(X,Y)\le K' \log(i(X,Y))+L'.
\end{displaymath}
\end{theo}

Given a basis $x$ of $F_n$ and a word $w\in F_n$, we denote by $|w|_x$ the length of the word $w$, when written in the basis $x$. Given two bases $x$ and $y$ of $F_n$, we denote by $|y|_x$ the maximal length of an element in $y$, when written in the basis $x$.

\begin{lemma} \label{lemma2.3}
Let $x=(x_i),y=(y_i),w=(w_i)$ be three bases of $F_n$, and suppose that there exists $v\in F_n$ such that for all $i\in[|1,n|]$, we have $x_i=vw_iv^{-1}$. Then $|v|_y\le |w|_y^2 |y|_x$.
\end{lemma}

\noindent\textit{Proof} : If $v$ can be written as a word $W(w_i)$, then we also have $v=W(x_i)$. Hence $|v|_x=|v|_w$. Besides, for all $i\in [|1,n|]$, we have $w_i=v^{-1}x_iv$, so there exists $i_0\in [|1,n|]$ such that $|w_{i_0}|_x=2|v|_x+1$. Hence

\begin{align*}
|v|_y&\le  |v|_w|w|_y\\
&= |v|_x|w|_y\\
&\le |w_{i_0}|_x|w|_y\\
&\le |w_{i_0}|_y|y|_x|w|_y\\
&\le |w|_y^2 |y|_x.
\end{align*}
\qed

Let $\Gamma_1,\Gamma_2\in CV_n$. In \cite[Lemma 3.4]{FM}, Francaviglia and Martino proved the existence of a difference of markings from $\Gamma_1$ to $\Gamma_2$ with minimal Lipschitz constant, as an easy application of the Arzelà-Ascoli theorem. We show, when $\Gamma_1$ and $\Gamma_2$ are roses, that we can choose the difference of markings to send the vertex of $\Gamma_1$ to the vertex of $\Gamma_2$, without changing too much its Lipschitz constant. 

\begin{lemma} \label{lemma2.4}
Let $\Gamma_1,\Gamma_2\in CV_n$ be two roses with all petals having length $\frac{1}{n}$. There exists a morphism $f:\Gamma_1\to\Gamma_2$ whose Lipschitz constant is no greater than $2\Lambda(\Gamma_1,\Gamma_2)$.
\end{lemma}

\noindent \textit{Proof} : Denote by $v_1$ (resp. $v_2$) the vertex of $\Gamma_1$ (resp. $\Gamma_2$). Let $g:\Gamma_1\to\Gamma_2$ be a difference of markings from $\Gamma_1$ to $\Gamma_2$ whose Lipschitz constant is equal to the stretching factor $\Lambda(\Gamma_1,\Gamma_2)$. Fix a path $\gamma\subset\Gamma_2$ from $v_2$ to $g(v_1)$ having minimal length in $\Gamma_2$. In particular, the path $\gamma$ has length no greater than $\frac{1}{2n}$. We define a morphism $f:\Gamma_1\to\Gamma_2$ by sending each petal $e$ of $\Gamma_1$ linearly to the path obtained by tightening the concatenation $\gamma g(e) \overline{\gamma}$. The image $f(e)$ has length at most $\frac{\Lambda(\Gamma_1,\Gamma_2)+1}{n}$, which is no more than $2\frac{\Lambda(\Gamma_1,\Gamma_2)}{n}$ as $\Lambda(\Gamma_1,\Gamma_2)\ge 1$. This implies that the Lipschitz constant of the morphism $f$ is no greater than $2\Lambda(\Gamma_1,\Gamma_2)$.
\qed
\\
\\
In the sequel, we will call a morphism given by lemma \ref{lemma2.4} \emph{quasi-optimal}.
\\
\\
We fix a standard basis $x=(x_1,\dots,x_n)$ of $F_n$. Let $\Gamma_1,\Gamma_2\in CV_n$ be two roses, the rose $\Gamma_1$ being the standard rose. There is a natural basis associated to any morphism $\sigma':\Gamma_2\to\Gamma_1$ : each petal of $\Gamma_2$ is labelled with the word of $F_n$ defined by its image by $\sigma'$. Conversely, a basis $(y_i)$ of $F_n$ defines a morphism from $\Gamma_2$ to $\Gamma_1$ by subdividing the petals into $|y_i|_x$ segments of length $\frac{1}{n|y_i|_x}$, and mapping them linearly to the corresponding petal in $\Gamma_1$.

\begin{lemma} \label{lemma2.5}
For all $n\ge 2$, there exists $C=C(n)\in\mathbb{R}$ such that the following holds.
\\
Let $\Gamma_1,\Gamma_2\in CV_n$ be two roses with petal lengths $\frac{1}{n}$. We assume that $\Gamma_1$ is the standard rose with petals labelled by $x_1,\dots,x_n$. Let $\sigma':\Gamma_2\to\Gamma_1$ be a quasi-optimal morphism, and let $y$ be the associated basis. Then $|y|_x\le 2\Lambda(\Gamma_2,\Gamma_1)$ and $|x|_y\le C\Lambda(\Gamma_1,\Gamma_2)^C$. 
\end{lemma}

\noindent\textit{Proof} : For all $i\in [|1,n|]$, the loop $y_i$ is subdivided into $|y_i|_x$ subsegments of length $\frac{1}{n|y_i|_x}$, and each of these subsegments is mapped by $\sigma'$ to a loop in $\Gamma_1$ of length $\frac{1}{n}$. As the Lipschitz constant of $\sigma'$ is no greater than $2\Lambda(\Gamma_2,\Gamma_1)$, we have $|y_i|_x\le 2\Lambda(\Gamma_2,\Gamma_1)$ for all $i\in [|1,n|]$, hence $|y|_x\le 2\Lambda(\Gamma_2,\Gamma_1)$.
\\
\\
Let $x'=(x'_i)$ be a basis associated to a quasi-optimal morphism from $\Gamma_1$ to $\Gamma_2$. Then there exists $v\in F_n$ such that for all $i\in [|1,n|]$, we have $x'_i=v^{-1}x_iv$. By the previous argument, we have $|y|_x\le 2\Lambda(\Gamma_2,\Gamma_1)$ and $|x'|_y\le 2\Lambda(\Gamma_1,\Gamma_2)$. In addition, we have $|x|_y\le 2 |v|_y+|x'|_y$, so lemma \ref{lemma2.3} implies that $|x|_y\le 2|x'|_y^2|y|_x+|x'|_y$. Using proposition \ref{prop1.1}, we thus have a polynomial bound on $|x|_y$ in terms of $\Lambda(\Gamma_1,\Gamma_2)$.
\qed
\\
\\
\indent We now relate the stretching factor from a point in outer space to two points that are close in outer space, and give a similar estimate for intersection numbers. This will allow us to deal only with the case when $X$ and $Y$ are roses in the proof of theorem \ref{theo2.2}. 

\begin{lemma} \label{lemma2.6}
For all $n\ge 2$ and $\epsilon>0$, there exists $C=C(n,\epsilon)\in\mathbb{R}$ such that if $A\in CV_n$ and if $B,B'\in CV_n$ are two points in the $\epsilon$-thick part of outer space whose simplices share the simplex of a rose as a face, then $\Lambda(A,B')\le C \Lambda(A,B)$ and $\Lambda(B',A)\le C\Lambda(B,A)$.
\end{lemma}

\noindent\textit{Proof} : Let $K_0$ denote the diameter of the $\epsilon$-thick part of the star of a rose simplex in outer space (i.e. the supremum of $d(A,B)$ for points $A$,$B$ lying inside it). This does not depend on the choice of the rose because $Out(F_n)$ acts transitively on the simplices of roses, and the metric is $Out(F_n)$-invariant. Besides, as the $\epsilon$-thick part of the star of a rose simplex is compact, we get that $K_0$ is finite. If $B,B'\in CV_n$ are two points in the $\epsilon$-thick part of outer space whose simplices share the simplex of a rose as a face, then $d(B,B')\le K_0$. By the triangle inequality, for all $A\in CV_n$, we have $d(A,B')\le d(A,B)+K_0$, hence $\Lambda(A,B')\le e^{K_0}\Lambda(A,B)$. A similar argument shows that $\Lambda(B',A)\le e^{K_0}\Lambda(B,A)$.
\qed

\begin{lemma} \label{lemma2.7}
For all $n\ge 2$, there exist $C=C(n)\in\mathbb{R}$ and $D=D(n)\in\mathbb{R}$ such that if $A\in CV_n$ and if $B,B'\in CV_n$ are two points in outer space whose simplices share the simplex of a rose as a face, then $i(A,B')\le Ci(A,B)+D$.
\end{lemma}

\noindent\textit{Proof} : This is a consequence of \cite[Lemma 2.7]{BBC} and the $Out(F_n)$-invariance of intersection numbers. Indeed, if the simplices of $B$ and $B'$ share the simplex of a rose as a face, then there exists $B_0,B_0'\in CV_n$ whose simplices share the simplex of the standard rose as a face, and $\Phi\in Out(F_n)$ such that $B=\Phi(B_0)$ and $B'=\Phi(B_0')$. Now notice that there are only finitely many simplices that have the simplex of the standard rose as a face. 

\indent We now give a geometric proof of this proposition in the model of spheres of outer space, that does not rely on \cite[Lemma 2.7]{BBC}. Let $B$ and $B'$ be two sphere systems whose simplices share the simplex of a rose as a face. We first assume that $B\subset B'$, and that $A$ is in normal form with respect to $B'$. It is easy to check that a graph with fundamental group $F_n$ has at least $n$ edges, and at most $3n-3$ edges. Hence the sphere system $B'$ can have at most $2n-3$ more spheres than $B$. Let $P$ be a component of $M_n-B$. As $A$ is in normal form with respect to $B'$, each component of $A\cap P$ intersects each sphere of $B'$ at most once. The number of components of $A\cap P$ is at most $i(A,B)+2n-3$. Therefore $B'$ has at most $(2n-3)(i(A,B)+2n-3)$ more intersection circles with $A$ than $B$, so $i(A,B')\le (2n-2)i(A,B)+(2n-3)^2$.

More generally, let $R$ be a rose such that the simplices of $B$ and $B'$ share the simplex of $R$ as a face. Then by the previous paragraph, we have $i(A,B)\le C i(A,R)$ and $i(A,R)\le i(A,B')$, so the result follows.
\qed
\\
\\
\indent Using the construction of Behrstock, Bestvina and Clay, and proposition \ref{prop1.2}, we prove the following estimate between intersection numbers and stretching factors. This is the left-hand side inequality of theorem \ref{theo2.2}. In the proof, we will use the notations introduced in section \ref{sec1.2}, in the paragraph preceding proposition \ref{prop1.2}. 

\begin{prop} \label{prop2.8}
Given $n\ge 2$ and $\epsilon>0$, there exists $C=C(n,\epsilon)>0$ such that for all points $T_1$ and $T_2$ in the $\epsilon$-thick part of $CV_n$, we have 
\begin{displaymath}
i(T_1,T_2)\le C\Lambda(T_1,T_2)^C,
\end{displaymath}
\noindent i.e.
\begin{displaymath}
\frac{1}{C}\log i(T_1,T_2)-\frac{\log C}{C}\le d(T_1,T_2).
\end{displaymath}
\end{prop}

\noindent \textit{Proof} : Let $\Gamma_1=T_1/F_n$ and $\Gamma_2=T_2/F_n$. Assume first that $\Gamma_1$ and $\Gamma_2$ are roses, all of whose petals have length $\frac{1}{n}$. Let $\sigma':\Gamma_2\to\Gamma_1$ be a quasi-optimal morphism, and label the petals of $\Gamma_2$ by the corresponding basis. By lemma \ref{lemma2.5}, after subdividing the petals of $\Gamma_2$ in at most $2\Lambda(\Gamma_2,\Gamma_1)$ segments, the morphism $\sigma'$ sends each segment of $\Gamma_2$ to a petal of $\Gamma_1$. Hence for every edge $e\subset\Gamma_1$, the set $\Sigma_e$ has cardinality at most $2n\Lambda(\Gamma_2,\Gamma_1)$. Let $e\subset\Gamma_1$ be an edge, and $q\in\Sigma_e$. To construct $\gamma_q$, start by joining $*_2$ to $q$ by a path $\gamma$, staying in one petal. The $\sigma'$-image of this path in $\Gamma_1$ crosses at most $2\Lambda(\Gamma_2,\Gamma_1)$ petals (counted with multiplicities) denoted by $e_1,\dots,e_k$ (the $e_i$ are not necessarily distinct) before crossing either $e_+$ or $e_-$. Write each $e_i$ as a word $e_i^1,\dots,e_i^{j_i}$ in the petals of $\Gamma_2$. Lemma \ref{lemma2.5} ensures that $j_i$ is polynomially bounded in $\Lambda(\Gamma_1,\Gamma_2)$. We form $\gamma_q$ by first crossing the petals $\overline{e_k^{j_k}},\dots,\overline{e_1^1}$, then crossing the path $\gamma$, and tightening. So the length of $\gamma_q$ is polynomially bounded in $\Lambda(T_1,T_2)$ and $\Lambda(T_2,T_1)$. Using proposition \ref{prop1.1}, we get that the length of $\gamma_q$ is bounded by a polynomial function of $\Lambda(T_1,T_2)$, and so is the number of edges in the tree $T_e$. Proposition \ref{prop1.2} then implies the proposition in the case when $\Gamma_1$ and $\Gamma_2$ are roses with all petals having length $\frac{1}{n}$. The general case follows from lemmas \ref{lemma2.6} and \ref{lemma2.7}, as one can get a rose from any graph by collapsing a maximal tree.
\qed
\\
\\
\indent We now use the geometric interpretation of intersection numbers to prove a converse estimate. Let $Y$ be a simple sphere system in $M_n$ dual to the standard rose, and $X$ be a simple sphere system dual to a rose in normal form with respect to $Y$. We identify each of the spheres in $Y$ with one of the generators of $F_n$. Let $G_X$ be a graph embedded in $M_n$ dual to the sphere system $X$. There is a natural basis associated to $G_X$ : each petal of $G_X$ is labelled by the word corresponding to the successive spheres in $Y$ it crosses.

\begin{prop} \label{prop2.9}
For all $n\ge 2$ and $\epsilon>0$, there exist $A,B\in\mathbb{R}$ such that for all points $X$ and $Y$ in the $\epsilon$-thick part of $CV_n$, we have
\begin{center}
$\frac{1}{A}\Lambda(X,Y)-B\le i(X,Y)$.
\end{center}
\noindent Hence for all $n\ge 2$ and $\epsilon>0$, there exists $B'\in\mathbb{R}$ such that for all points $X$ and $Y$ in the $\epsilon$-thick part of $CV_n$, we have
\begin{center}
$d(X,Y)\le\log i(X,Y)+B'$.
\end{center}
\end{prop}

\noindent \textit{Proof} : Let $X,Y\in CV_n$. We first assume that $X$ and $Y$ correspond to roses whose petals have length $\frac{1}{n}$, the rose corresponding to $Y$ being the standard rose of $F_n$. We claim that we can find a graph $G_X$ embedded in $M_n$ and dual to $X$ which crosses at most twice each connected component of $Y-X$. Indeed, let $G'_X$ be a graph embedded in $M_n$ and dual to $X$. We first homotope $G'_X$ so that it does not cross any of the intersection circles in $X\cap Y$. Suppose that an edge $e$ of $G'_X$ crosses a component of $Y-X$ in three points $x_1,x_2,x_3$, and let $\gamma$ (resp. $\gamma'$) be the subpath of $e$ from $x_1$ to $x_2$ (resp. $x_2$ to $x_3$). As $e$ crosses $X$ exactly once, one of the paths $\gamma$ and $\gamma'$ does not cross $X$. Without loss of generality, we assume that $\gamma$ does not cross $X$. As $X$ is a simple sphere system, the path $\gamma$ stays in a simply connected region of $M_n$, so we can homotope $\gamma$ to a path from $x_1$ to $x_2$ that remains in one component of $Y-X$. We can then slightly homotope the new edge $e$ so that it crosses this component of $Y-X$ at most twice. 

Denote by $C_0$ the dimension of $CV_n$, i.e. the maximal number of spheres in a simple sphere system. Let $s'$ be a sphere in $Y$. By induction on the number of intersection circles between $X$ and $s'$, one gets that the number of connected components of $s'-X$ is equal to $i(X,s')+1$. Hence each of the edges in $G_X$ crosses $Y$ at most $2(i(X,Y)+C_0)$ times. Denote by $x'$ the basis of $F_n$ associated to $G_X$, then we have $|x'|_x\le 2(i(X,Y)+C_0)$. The morphism from $G_X$ to $G_Y$ defined by the basis $x'$ has Lipschitz constant $|x'|_x$, so $\Lambda(X,Y)\le |x'|_x$. The claim follows in the case when $X$ and $Y$ correspond to roses with all petals having length $\frac{1}{n}$, the rose corresponding to $Y$ being the standard rose of $F_n$. The general case follows from lemmas \ref{lemma2.6} and \ref{lemma2.7}.

So there exist $A,B\in\mathbb{R}$ such that for all $X,Y$ in the $\epsilon$-thick part of outer space, we have 
\begin{center}
$\Lambda(X,Y)\le A i(X,Y)+AB$.
\end{center}
\noindent If $i(X,Y)\ge 1$, we thus have
\begin{center}
$\Lambda(X,Y)\le (A+AB)i(X,Y)$,
\end{center}
\noindent i.e.
\begin{center}
$d(X,Y)\le \log i(X,Y) + \log(A+AB)$.
\end{center}
\noindent If $i(X,Y)=0$, then $X$ and $Y$ are compatible, so $d(X,Y)\le K_0$, where $K_0$ is the maximal diameter of the $\epsilon$-thick part of a closed simplex in outer space. So letting $B'=\max(\log(A+AB),K_0)$ gives the last inequality in the proposition.
\qed
\\
\\
Theorem \ref{theo2.2} follows from propositions \ref{prop2.8} and \ref{prop2.9}.

\section{Metric properties of the combing path} \label{sec3}

Given $A,B\in CV_n$, the combing path $\gamma$ from $B$ to $A$ is a piecewise linear path, of the form $B=A_N,\dots,A_1=A$. We recall from the introduction that we denote by $l(\gamma)$ the \emph{length} of $\gamma$, defined to be 

\begin{displaymath}
l(\gamma):=\sum_{i=1}^{N-1} d(A_i,A_{i+1}).
\end{displaymath}

\noindent The goal of this section is to prove the main theorem of this paper, which states that combing paths make definite progress in outer space, provided they remain in the $\epsilon$-thick part for some $\epsilon>0$.

\begin{theo} \label{theo3.1}
\textit{For all $n\ge 2$ and $\epsilon>0$, there exist $K,L\in\mathbb{R}$ such that the following holds.}

\noindent \textit{Let $A,B\in CV_n$ be such that the combing path $\gamma$ from $B$ to $A$ remains in the $\epsilon$-thick part of outer space. Then}

\begin{center}
\textit{$\frac{l(\gamma)}{K}-L\le d(A,B)\le l(\gamma).$}
\end{center}
\end{theo}

\subsection{A few facts about combing paths} \label{sec3.1}

\begin{figure} 
\begin{center}
\includegraphics[width=12cm,height=3cm]{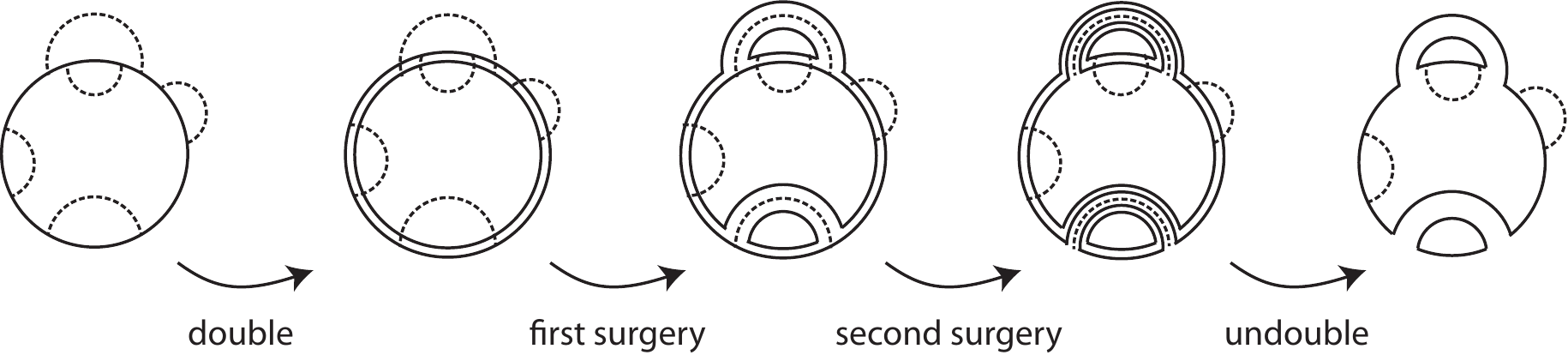}
\caption{Distribution of intersection circles during a double surgery step}
\label{fig3}
\end{center}
\end{figure}

\begin{figure} 
\begin{center}
\includegraphics[width=9cm,height=3cm]{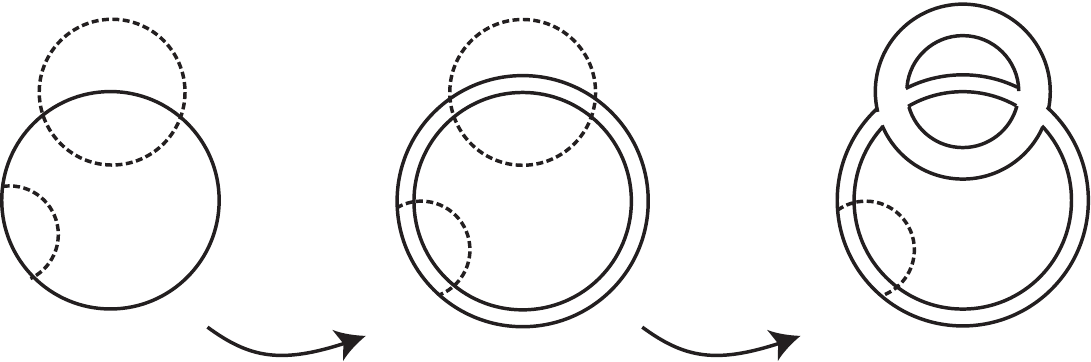}
\caption{Case of a surgery along the last intersection circle with a sphere in $A$}
\label{fig4}
\end{center}
\end{figure}

\noindent We start by collecting some facts about combing paths, which follow from the construction of these paths described in section \ref{sec1.4}. Let $A,B\in CV_n$, and let $A=A_0,\dots,A_N=B$ be the vertices in the combing path from $B$ to $A$. In this setting, a double surgery step as described in section \ref{sec1.4} consists of passing from $A_{N-i}$ to $A_{N-i-2}$ for some even $i\in[|0,N|]$, the sphere system $A_{N-i-1}$ being the sphere system obtained from $A_{N-i}$ by doubling, performing one surgery, and projecting to outer space. 
\\
\\
\textit{Fact 1} : For all $i\in[|0,N|]$ such that $N-i$ is even, the combing path from $A_i$ to $A$ is a subpath of the combing path from $B$ to $A$. 
\\
\\
We want to understand the evolution of intersection numbers along combing paths. We describe the evolution of intersection circles between $A$ and $A_i$ for $i\in [|0,N|]$. On Figure \ref{fig3}, we draw in dodded lines the intersection circles between $A$ and a sphere system $B$, and we look at how they evolve when performing a double surgery step. Note that doubling the sphere system $B$ may (at most) double the number of intersection circles with $A$, but undoubling the sphere system at the end ensures that intersection cirles get distributed over the created spheres. However, as we see in Figure \ref{fig4}, this distribution does not occur when we are performing surgery along the last intersection circle between a sphere in $A$ and $B$. We collect these observations in the following two facts. 
\\

\begin{figure} 
\begin{center}
\includegraphics[width=8cm,height=3cm]{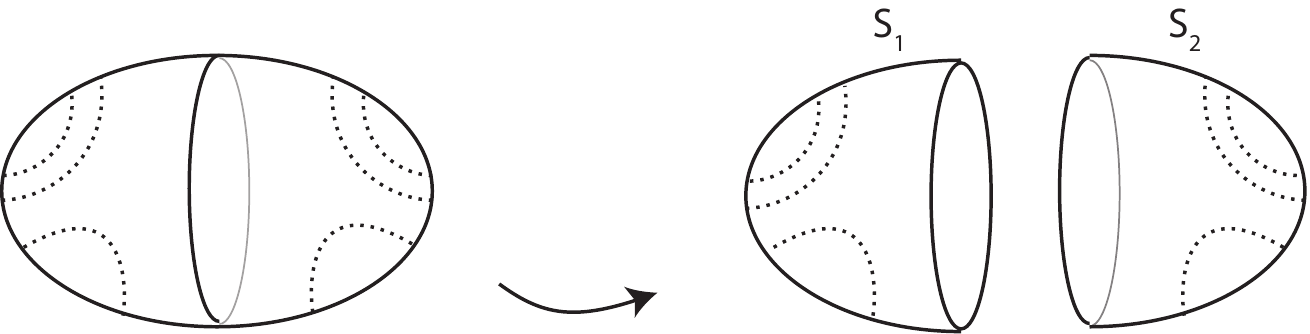}
\caption{Understanding the case when only one sphere is created}
\label{fig5}
\end{center}
\end{figure}

\textit{Fact 2} : After performing a double surgery step on a sphere $s\in B$, the intersection circles in $s\cap A$ are distributed over the spheres that come from these surgeries on $s$ (some are even deleted), except possibly when performing surgery along the last intersection circle of a sphere in $A$. However, note that this exceptional case cannot occur more than $C_0$ times, where $C_0$ is the dimension of $CV_n$, i.e. the maximal number of spheres in a sphere system in $CV_n$. In particular, this implies that for $k\in[|0,\lfloor\frac{N}{2}\rfloor|]$, we have $i(A,A_{N-2(k+1)})\le i(A,A_{N-2k})$, except for at most $C_0$ values of $k$.
\\
\\
\textit{Fact 3} : For all $k\in [|0,N-1|]$, we have $i(A,A_k)<2i(A,A_{k+1})$ (the inequality is strict because the intersection circle used to perform surgery is removed).
\\
\\
Finally, we will have to understand what happens when only one sphere is created when performing surgery. Suppose that after performing surgery as on figure \ref{fig5}, we get only one sphere. This means that either one of the spheres $S_1$ or $S_2$ is trivial, or that both get identified. The first case is impossible because it would contradict the fact that $B$ is in normal form with respect to $A$, so the spheres $S_1$ and $S_2$ are parallel. This implies in particular (as on figure \ref{fig5}) that the pattern of intersection circles between $A$ and each of these spheres is the same. So the spheres $S_1$ and $S_2$ have at least twice fewer intersection circles with $A$ than $B$ had. More generally, we get the following fact.
\\

\textit{Fact 4} : As a consequence of normal form, it is impossible that after performing a double surgery step, all the spheres you get except one are trivial. So you cannot get only one sphere, except if (at least) two of the spheres coming from the initial sphere have been identified. In particular, they have at least twice fewer intersections with $A$ than the initial sphere did.

\subsection{Growth of intersection numbers along combing paths and end of the proof} \label{sec3.2}

In order to prove theorem \ref{theo3.1}, we first determine the growth of intersection numbers between vertices in the combing path. The following proposition may be seen as an analog for combing paths of the result of Behrstock, Bestvina and Clay establishing exponential growth of intersection numbers along axes of fully irreducible automorphisms of $F_n$ in outer space (\cite{BBC}). In addition, Behrstock, Bestvina and Clay give the exact growth rate of intersection numbers, in terms of the Perron-Frobenius eigenvalues of the fully irreducible automorphism and its inverse. Note that combining the fact that folding paths are geodesics in $CV_n$ \cite[Theorem 5.5]{FM} with proposition \ref{prop2.9} shows that intersection numbers also grow exponentially along folding paths that remain in the $\epsilon$-thick part of outer space for some $\epsilon>0$.

\begin{prop} \label{prop3.2}
\textit{For all $n\ge 2$ and $\epsilon>0$, there exist $C_1, C_2\in (1,+\infty)$ such that the following holds.
\\ Let $A,B\in CV_n$ be such that the combing path from $B$ to $A$ stays in the $\epsilon$-thick part of $CV_n$, and let $A=A_0,\dots,A_N=B$ be the vertices of this path. Assume that $N\ge 3$, then }

\begin{displaymath} 
C_1^{N}\le i(A,B)\le C_2^{N}.
\end{displaymath}
\end{prop}

\noindent \textit{Proof of the upper bound} : The combing path from $B$ to $A$ is a piecewise linear path, each piece staying in one closed simplex ; we denote by $A=A_0,\dots,A_N=B$ its vertices. Let $C,D$ be the constants given by lemma \ref{lemma2.7}. We prove by induction on $k\in [|0,N|]$ that $i(A,A_k)\le (C+D)^{k}$. This is obviously true for $k=0$. Assume that $i(A,A_k)\le (C+D)^{k}$. If $i(A,A_k)=0$, then by lemma \ref{lemma2.7}, we have $i(A,A_{k+1})\le D\le (C+D)^{k+1}$. If $i(A,A_k)>0$, then by lemma \ref{lemma2.7}, we have $i(A,A_{k+1})\le (C+D) i(A,A_{k})$, so by the induction hypothesis we have $i(A,A_{k+1})\le (C+D)^{k+1}$.
\qed
\\
\\
We now prove the lower bound in proposition \ref{prop3.2}. Let $0\le i\le j\le N$ be such that $N-i$ and $N-j$ are even, and such that when performing the surgeries leading from $A_j$ to $A_i$, you never perform surgery along the last intersection circle of a sphere in $A$. Let $M=\frac{j-i}{2}$. For $p\in [|0,M|]$, let $S_p^1,\dots,S_p^{k_p}$ be the spheres in $A_{j-2p}$, i.e. the spheres you get from $A_{j}$ after performing $2p$ successive surgery steps (i.e. $p$ double surgery steps). 
\\
\\
We define a function $\theta$ on finite tuples of integers by 

\begin{displaymath}
\theta (i_1,\dots,i_k)=
\left\lbrace
\begin{array}{ccc}
i_1+1  & \text{\ if \ } k>1 \text{\ and \ } i_1=i_2=\dots =i_k\\
\max\{i_1,\dots ,i_k\} & \text{\ otherwise \ } 
\end{array}\right.
\end{displaymath}

\noindent For all $p\in[|0,M|]$ and all $q\in [|1,k_p|]$, we associate to the sphere $S_p^q$ an integer $X(S_p^q)$, by downward induction on $p$, in the following way. This integer will help us count the intersection circles between $S_p^q$ and $A_i$.
\\
\\
$\bullet$ \ For all $l\in [|1,k_{M}|]$, let $X(S_{M}^l)=0$.
\\
\\
$\bullet$ \ Let $p\in [|1,M|]$ and $q\in [|1,k_{p-1}|]$. Assume that we have defined $X(S_p^l)$ for all $l\in [|1,k_p|]$. We want to define $X(S_{p-1}^q)$.
\\
\\
\textit{Case 1} : When performing surgery on the sphere system $A_{j-2(p-1)}$ along $A$, no surgery is performed on $S_{p-1}^q$. Then $S_{p-1}^q=S_{p}^l$ for some $l\in [|1,k_p|]$, and we let $X(S_{p-1}^q)=X(S_p^l)$.
\\
\\
\textit{Case 2} : The sphere $S_{p-1}^q$ intersects the sphere system $A$ in at least one circle that bounds an innermost disk in $A$. In addition, after performing a double surgery step on the sphere $S_{p-1}^q$, only one nontrivial sphere $S_p^{\alpha}$ is created (after identifying parallel spheres). Then we let $X(S_{p-1}^q)=X(S_p^{\alpha})+1$. 
\\
\\
\textit{Case 3} : The sphere $S_{p-1}^q$ intersects the sphere system $A$ in at least one circle that bounds an innermost disk in $A$. In addition, after performing a double surgery step on the sphere $S_{p-1}^q$, at least two nontrivial spheres are created (after identifying parallel spheres). We denote by $S_{p}^{\alpha_1},\dots,S_{p}^{\alpha_l}$ the created spheres. Then we let $X(S_{p-1}^q)=\theta (X(S_{p}^{\alpha_1}),\dots, X(S_{p}^{\alpha_l}))$.
\\
\\
\indent Define a sequence $u$ by $u_0=0$ and $u_n=2^{n-1}$ for all $n\ge 1$. The following lemma gives a way of comparing the growth of intersection numbers along the combing path with the exponential growth of $u$.

\begin{lemma} \label{lemma3.3}
Let $p\in [|0,M|]$, and $q\in [|1,k_p|]$. Then the number of intersection circles between $S_p^q$ and the sphere system $A$ is at least $u_{X(S_p^q)}$.
\end{lemma}  

\noindent \textit{Proof} : The proof is by downward induction on $p$. The result is obvious when $p=M$. Now assume that for all $l\in [|1,k_p|]$, the number of intersection circles between $S_p^l$ and the sphere system $A$ is at least $u_{X(S_p^l)}$. We want to show that for all $q\in [|1,k_{p-1}|]$, the number of intersection circles between $S_{p-1}^q$ and $A$ is at least $u_{X(S_{p-1}^q)}$.

In case 1, let $l\in [|1,k_p|]$ be such that $S_{p-1}^q=S_{p}^l$. By the induction hypothesis, the number of intersection circles between $S_{p}^l$ and $A$ is greater than $u_{X(S_p^l)}$. But by definition, we have $X(S_{p-1}^q)=X(S_p^l)$, so the number of intersection circles between $S_{p-1}^q=S_{p}^l$ and $A$ is greater than $u_{X(S_{p-1}^q)}$.

Suppose now that we are in case 2, and let $S_p^{\alpha}$ be the unique sphere that is created from $S_{p-1}^q$ after performing a double surgery step. By fact 4 of the previous section, the number of intersection circles between $S_{p-1}^q$ and $A$ is at least twice the number of intersection circles between $S_p^{\alpha}$ and $A$. By the induction hypothesis, the number of intersection circles between $S_p^{\alpha}$ and $A$ is at least $u_{X(S_p^{\alpha})}$. So the number of intersection circles between $S_{p-1}^q$ and $A$ is at least $2u_{X(S_p^{\alpha})}$, which is greater than $u_{X(S_p^{\alpha})+1}=u_{X(S_{p-1}^q)}$ if $X(S_{p}^{\alpha})>0$. If $X(S_{p}^{\alpha})=0$, then $X(S_{p-1}^q)=1$, and the number of intersection circles between $S_{p-1}^q$ and $A$ is at least $1$ by the first assumption made in case 2.

Finally, suppose that we are in case 3, and let $S_{p}^{\alpha_1},\dots,S_{p}^{\alpha_l}$ be the spheres created from $S_{p-1}^q$ after performing a double surgery step. By fact 2 of the previous section, the number of intersection circles between $S_{p-1}^q$ and $A$ is greater than the sum of the number of intersection circles between $A$ and the spheres $S_p^{\alpha_j}$ for $j\in [|1,l|]$. By the induction hypothesis, this sum is at least equal to $u_{X(S_{p}^{\alpha_1})}+\dots +u_{X(S_{p}^{\alpha_l})}$. This is greater than $u_{\max_j X(S_{p}^{\alpha_j})}$, and if all $u_{X(S_{p}^{\alpha_j})}$ are equal and different from $0$, it is greater than $u_{X(S_{p}^{\alpha_1})+1}$. In both cases, the sum is greater than $u_{\theta (X(S_{p}^{\alpha_1}),\dots, X(S_{p}^{\alpha_l}))}=u_{X(S_{p-1}^q)}$. So the number of intersection circles between $S_{p-1}^q$ and $A$ is at least $u_{X(S_{p-1}^q)}$. If all $X(S_{p}^{\alpha_j})$ are equal to $0$, then $u_{X(S_{p-1}^q)}=u_{\theta (X(S_{p}^{\alpha_1}),\dots, X(S_{p}^{\alpha_l}))}=1$, and the number of intersection circles between $S_{p-1}^q$ and $A$ is at least $1$ by the first assumption made in case 3.
\hfill $\square$
\\
\\
\indent As the sequence $u$ grows exponentially, our goal is now to prove that the integers $X(S_p^q)$ grow linearly along the combing path. We denote by $w(s)$ the weight of a sphere $s$. For $p\in [|0,M|]$, we define
\begin{displaymath}
N(p)=\sum_{q=1}^{k_p} w(S_p^q)X(S_p^q).
\end{displaymath}

\begin{lemma} \label{lemma3.4}
There exists $C_3>0$, such that for all $p\in [|0,M-1|]$, we have
\begin{displaymath}
N(p)-N(p+1)\ge\frac{1}{C_3}\max_s(w(s)),
\end{displaymath}
\noindent the maximum being taken over all spheres $s\in A_{j-2p}$ that get nontrivially subdivided when performing surgery along $A$ on the sphere system $A_{j-2p}$. 
\end{lemma}

Recall that in the definition of the combing path, when a sphere $s$ gets subdivided, its weight is transferred equally among all the nontrivial spheres that come from it. If a sphere $s'$ is obtained from $s$ after a double surgery step, we denote by $w_s(s')$ the part of the weight of $s'$ that comes from the sphere $s$.
\\
\\
\noindent \textit{Proof of lemma \ref{lemma3.4}} : Let $C_0$ denote the dimension of $CV_n$, i.e. the maximal number of edges in a graph in $CV_n$. Let $s$ be a sphere that gets nontrivially subdivided when performing a double surgery step along $A$ on the sphere system $A_{j-2p}$, and let $S_{p+1}^{\alpha_1},\dots,S_{p+1}^{\alpha_l}$ be the spheres obtained from $s$ after a double surgery step. These spheres get a weight at least equal to $\frac{w(s)}{C_0}$ from $s$. 
\\
\\
\noindent $\bullet$ \ If there is only one such sphere, then $X(s)=X(S^{\alpha_1}_{p+1})+1$ (case 2 of the definition of the integers $X$), so

\begin{center}
$w(s)X(s)-w_s(S_{p+1}^{\alpha_1})X(S_{p+1}^{\alpha_1})=w(s)$.
\end{center}

\noindent $\bullet$ \ If there are at least two such spheres, and if $X(S_{p+1}^{\alpha_1})=\dots=X(S_{p+1}^{\alpha_l})$, then by case 3 of the definition of the integers $X$ and the definition of $\theta$, we have $X(s)=X(S^{\alpha_1}_{p+1})+1$, so

\begin{displaymath}
w(s)X(s)-\sum_{i=1}^l w_s(S_{p+1}^{\alpha_i})X(S_{p+1}^{\alpha_i})=w(s).
\end{displaymath}

\noindent $\bullet$ \ Finally, if there are at least two such spheres, and if there exist $i_0,i_1\in [|1,l|]$ with $X(S_{p+1}^{\alpha_{i_0}})<X(S_{p+1}^{\alpha_{i_1}})$, then by case 3 of the definition of the integers $X$ and the definition of $\theta$, we have $X(s)=\max_i\{X(S^{\alpha_i}_{p+1})\}$. Given a sphere $s\in A_{j-2p}$, define

\begin{displaymath}
N_s=\sum_{i=1}^l w_s(S_{p+1}^{\alpha_i})X(S_{p+1}^{\alpha_i}).
\end{displaymath}

\noindent We have

\begin{displaymath}
\begin{array}{rl}
w(s)X(s)-N_s &\ge w(s)X(s)-\sum_{i\ne i_0} w_s(S_{p+1}^{\alpha_i})X(s)-w_s(S_{p+1}^{\alpha_{i_0}})X(S_{p+1}^{\alpha_{i_0}}) \\
&\ge w(s)X(s)-\sum_{i\ne i_0} w_s(S_{p+1}^{\alpha_i})X(s)-w_s(S_{p+1}^{\alpha_{i_0}})(X(s)-1)\\
& = w_s(S_{p+1}^{\alpha_{i_0}})\\
&\ge\frac{w(s)}{C_0}.
\end{array}
\end{displaymath}

\noindent In all cases, we get that 

\begin{displaymath}
w(s)X(s)-\sum_{i=1}^l w_s(S_{p+1}^{\alpha_i})X(S_{p+1}^{\alpha_i})\ge\frac{w(s)}{C_0}.
\end{displaymath}

\noindent Summing the previous inequality over all the spheres in $A_{j-2p}$ that get subdivided, we get that

\begin{displaymath}
\begin{array}{rl}
N(p)-N(p+1)&=\sum_s (w(s)X(s)-N_s)\\ 
&\geq \sum_s \frac{w(s)}{C_0}\\
&     \geq \frac{1}{C_0}\max_s w(s).
\end{array}
\end{displaymath}
\qed

\begin{lemma} \label{lemma3.6}
Assume that the combing path from $B$ to $A$ remains in the $\epsilon$-thick part of outer space. Then there exist $C_4,C_5\in\mathbb{R}$ such that during a sequence of $C_4$ consecutive surgeries, at least one sphere with weight greater than $\frac{\epsilon}{C_5}$ gets subdivided.
\end{lemma}

\noindent \textit{Proof} : Let $\Sigma$ be a sphere system which is a vertex of the combing path. For $C_5\in\mathbb{R}$ big enough, the set of spheres $\Sigma'\subset\Sigma$ having weights less than $\frac{\epsilon}{C_5}$ corresponds to a forest in the corresponding graph (otherwise, as the number of edges is bounded, there is a loop of length less than $\epsilon$, which contradicts the assumption that the combing path stays in the $\epsilon$-thick part of $CV_n$). So the sphere system $\Sigma-\Sigma'$ is simple. Performing surgery on spheres in $\Sigma'$ creates a new sphere system of the form $\Sigma-\Sigma'\cup\Sigma''$, which is compatible with $\Sigma-\Sigma'$, and in which all spheres having weights less than $\frac{\epsilon}{C_5}$ belong to $\Sigma''$. Let $C_4$ be the maximal number of simplices in $CV_n$ corresponding to sphere systems at distance at most $2$ in the spine of outer space, which is finite because the action of $Out(F_n)$ on the spine is cocompact. Suppose that we perform $C_4+1$ consecutive surgeries only on spheres of weights less than $\frac{\epsilon}{C_5}$. Then all the $C_4+1$ sphere systems we get are compatible with $\Sigma-\Sigma'$, so we get back to a simplex we had already visited. This is impossible since the combing path must end at $A$.
\qed
\\
\\
\textit{Proof of proposition \ref{prop3.2}} : Let $0\le i\le j\le N$ be such that $N-i$ and $N-j$ are even, and that you never perform surgery along the last intersection circle of a sphere in $A$ between $A_j$ and $A_i$. Combining lemmas \ref{lemma3.4} and \ref{lemma3.6}, we get the existence of $C>0$ depending only on $n$ and $\epsilon$ such that for all $p\in [|0,M|]$, we have

\begin{center}
$N(p)-N(p+C)\ge 1.$
\end{center}

\noindent By induction on $p$, this implies that for all $p\in[|0,\lfloor\frac{M}{C}\rfloor|]$, we have

\begin{displaymath}
\sum_{q=1}^{k_{M-pC}} w(S_{M-pC}^q)X(S_{M-pC}^q)\ge p. 
\end{displaymath}

\noindent In particular, as the weights of a sphere system sum to $1$, one of the numbers $X(S_{M-pC}^q)$ is at least $p$. Let $p=\lfloor\frac{M}{C}\rfloor$. Denote by $C_0$ the dimension of outer space. By facts 2 and 3 of the previous section, we have $i(A,B)\ge \frac{1}{2^{C_0}}i(A,A_{j-2(M-pC)})$, so lemma \ref{lemma3.3} ensures that $i(A,B)\ge\frac{1}{2^{C_0}} u_{\lfloor\frac{M}{C}\rfloor}$. Let $N_0:=2\lceil C\rceil$. If $j-i\ge N_0$ (i.e. $M\ge\lceil C\rceil$), then 

\begin{center}
$\begin{array}{rl}
i(A,B)&\ge 2^{\lfloor\frac{M}{C}\rfloor -1-C_0}\\
&=2^{\lfloor\frac{j-i}{2C}\rfloor -1-C_0}\\
&\ge 2^{-2-C_0} (2^{\frac{1}{2C}})^{j-i}
\end{array}$
\end{center}

Let $N_0':=(C_0+1)(N_0+2)$. If $N\ge N_0'$, then we can subdivide the combing path from $B$ to $A$ into at most $C_0+1$ pieces, in which no surgery occurs on the last intersection circle with a sphere in $A$. One of these pieces contains at least $\frac{N}{C_0+1}\ge N_0+2$ vertices. So we can find $0\le i\le j\le N$ satisfying the above condition, and such that $j-i\ge\frac{N}{C_0+1}-2\ge N_0$. Hence we get

\begin{center}
$i(A,B)\ge 2^{-2-C_0-\frac{1}{C}} (2^{\frac{1}{2C(C_0+1)}})^N$.
\end{center}

We can thus find $N_0''\in\mathbb{N}$ and $C'_1>1$ such that if $N\ge N_0''$, then we have $i(A,B)\ge C_1'^{N}$. If $3\le N\le N_0''$, then $i(A,B)\ge 2$ (otherwise after one single surgery, we would get a sphere system compatible with $\Sigma$), so letting $C_1:=\min(2^{\frac{1}{N_0''}},C'_1)$ gives the result.
\qed
\\
\\
\textit{Proof of theorem \ref{theo3.1}} : The right-hand side inequality is an obvious application of the triangle inequality.

Let $K'$,$L'$ be the constants given by theorem \ref{theo2.2}, and let $C_1$ be given by proposition \ref{prop3.2}. Denote by $K_0$ the diameter of the $\epsilon$-thick part of the star of a rose simplex in $CV_n$. Assume that the combing path $\gamma$ from $B$ to $A$ remains in the $\epsilon$-thick part of $CV_n$, and denote by $A=A_0,\dots, A_N=B$ its vertices. As two consecutive vertices lie in the closure of the star of a common rose, we have $l(\gamma)\le K_0 N$. 
 If $N\ge 3$, then by proposition \ref{prop3.2} we have 

\begin{center}
$\log i(A,B)\ge N \log C_1$,
\end{center} 

\noindent hence by theorem \ref{theo2.2} we have

\begin{center}
$d(A,B)\ge \frac{\log C_1}{K'}N-L'$,
\end{center}

\noindent which implies that

\begin{center}
$d(A,B)\ge \frac{\log C_1}{K_0K'}l(\gamma)-L'$.
\end{center}

If $N\le 2$, then the length of the combing path is bounded above by $2K_0$. We conclude by letting $K:=\frac{\log C_1}{K_0K'}$ and $L=\max(L',\frac{2\log C_1}{K'})$.
\qed

\bibliographystyle{amsplain}
\bibliography{Biblio}

\end{document}